\documentclass{article}[12pt]
\usepackage{amsfonts}
\usepackage{amsmath}
\usepackage{amsthm}
\usepackage{enumerate}
\usepackage[british]{babel}
\usepackage[T1]{fontenc}

\usepackage{array}

\newcommand{\cA}{\mathcal A}

\newcommand{\cP}{\mathcal P}
\newcommand{\cB}{\mathcal B}
\newcommand{\cL}{\mathcal L}
\newcommand{\cM}{\mathcal M}
\newcommand{\GF}{\mathrm{GF}}
\newcommand{\PG}{\mathrm{PG}}
\newcommand{\cC}{\mathcal C}

\theoremstyle{theorem}
\newtheorem{lemma}{Lemma}
\newtheorem{theorem}[lemma]{Theorem}
\theoremstyle{remark}
\newtheorem{remark}{Remark}

\begin{document}
\title{On the non--existence of certain hyperovals in dual Andr\'e  planes
  of  order $2^{2h}$}
\author{A. Aguglia and L. Giuzzi\thanks{Research supported by  the Italian
    Ministry MIUR, Strutture geometriche, combinatoria e loro
    applicazioni.}}
\date{}

\maketitle

\begin{abstract}
No  regular hyperoval
of the Desarguesian  affine plane $AG(2,2^{2h})$,  with $h>1$,
is inherited by a dual Andr\'e plane of order $2^{2h}$ with
dimension $2$ over its centre.
\end{abstract}
\section{Introduction}

The general question on existence of ovals in finite
non--Desarguesian planes is still open and appears to be difficult.
It has been shown by computer search that there exist some planes
of order $16$ without ovals; see \cite{PRS}.
On the other hand, ovals
have been constructed in several finite planes;
one of the most fruitful approaches in this search has been that of
inherited oval, due to Korchm\'aros \cite{kg1974,kg1976}.

Korchm\'aros' idea relies on the fact that any two planes $\pi_1$
and $\pi_2$ of the same order have the same number of points and
lines; thus  their point sets, as well as some lines,  may be
identified. If $\Omega$ is an oval of $\pi_1$, it might happen that
$\Omega$, regarded as a point set,
turns also out to be an oval of $\pi_2$,
although $\pi_1$ and $\pi_2$ differ in
some (in general several) point--line incidences; in
this case $\Omega$ is called
an \emph{inherited oval} of $\pi_2$ from $\pi_1$;
see also \cite[Page 728]{HCD}.

In practice, it is usually convenient to take
$\pi_1$  to be the Desarguesian affine plane $AG(2,q)$ of
order a prime power $q$.
The case in which $\pi_2$ is the Hall plane $H(q^2)$ of order
$q^2$ was investigated in \cite{kg1974}, and inherited ovals
were found. For $q$ odd, this also proves the existence of inherited
ovals in the dual plane of $H(q^2)$,  which is a Moulton plane
$M(q^2)$ of the same order.
\par
Moulton planes have been originally introduced in \cite{Mo}, by
altering some of the lines of a Desarguesian plane constructed
over the reals, while keeping the original point set fixed.
In particular, each line of the Moulton plane turns out to be
either a line of the original plane or the union of two half--lines
of different slope with one point in common.
\par
This construction, when considering planes of finite order $q^2$, may
be carried out as follows.
Let $||\cdot||$ denote  the norm
function
\[ ||\cdot||:\begin{cases}
  \GF(q^2)\rightarrow\GF(q) \\
  x \mapsto x^{q+1}
  \end{cases}
\]
Take a proper subset $U$ of $\GF(q)^{\star}$ and
consider the following operation defined over the set $\GF(q^2)$
\[
 a\odot b=\begin{cases}
  ab  & \mbox{if $||b||\not\in U$} \\
  a^qb & \mbox{if $||b||\in U$}.
\end{cases}
\]
The set $(\GF(q^2),+,\odot)$ is a pre-quasifield  which is a
quasifield for $U\neq \{1\}$. Every pre-quasifield coordinatizes a
translation plane; see \cite[Section 5.6]{JJB}.
In our case this translation plane is an  affine Andr\'e plane
$A(q^2)$ of order $q^2$ with dimension $2$ over its centre.
In the case in which $U$ consists of a single element of $GF(q^2)$ the
translation plane is the affine Hall plane of order $q^2$  and its
dual plane is the affine Moulton plane of order $q^2$.
For details on these planes see \cite{D,HP}.

Write $\cM_U(q^2)=(\cP,\cL)$ for the incidence
structure whose
the point--set $\cP$  is the same as that of
$AG(2,q^2)$, and whose  lines in $\cL$ are either of the form
\[[c]=\{P(x,y): x=c, y\in\GF(q^2)\}\] or \[[m,n]=\{P(x,y): y=m\odot
x+n\}.\]
The affine plane $\cM_U(q^2)$ is the dual of an affine Andr\'e plane
$A(q^2)$ of order $q^2$.  Completing $\cM_U(q^2)$ with its points at
infinity in the usual way gives a projective plane
$\overline{\cM_U(q^2)}$ called \emph{the projective closure of
$\cM_U(q^2)$}.

Write $\Phi=\{P(x,y): ||x||\not\in U\}$ and
 $\Psi=\{P(x,y): ||x||\in U \}$.
Clearly, $\cP=\Phi\cup\Psi$.

If  an arc $\cA$   of $PG(2,q^2)$ is in turn an arc in
$\overline{\cM_U(q^2)}$ then, $\cA$ is
\emph{an inherited arc} of $\overline{\cM_U(q^2)}$.

Any hyperoval of the Desarguesian projective plane $PG(2,q^2)$
 obtained from a conic by adding its nucleus is called
\emph{regular}.
Let  consider a  regular hyperoval $\Omega$ in $AG(2,q^2)$.
If $\Omega \subseteq \Phi$,
that is for each point $P(x,y)\in \Omega$ the norm of $x$ is an element
of $\GF(q)\setminus U$, then $\Omega$ is clearly
an inherited hyperoval of $\cM_U(q^2)$.

In \cite[Theorem 1.1]{AL}, it is proven that for $q>5$ an odd prime
power, any arc of the Moulton plane $\cM_t(q^2)$ with $t\in GF(q)$, obtained as
$\cC^{\star}=\cC\cap\Phi$, where $\cC$ is an ellipse in $AG(2,q^2)$ is complete.

In this paper the  case where $q$ is even and $|U| < \frac{q}{4}-1$ is
addressed.   We prove the following.

\begin{theorem} \label{conf}
Suppose $\Omega$ to be a regular hyperoval of $AG(2,2^{2h})$, with $h>
1$.
 Then, $\Omega^*=\Omega \cap \Phi$ is a complete arc in the
 projective closure  of $\cM_U(2^{2h})$.
\end{theorem}
\begin{theorem}
\label{maint}
No regular hyperoval $\Omega$ of  $AG(2,2^{2h})$, with $h>1$,
is an inherited hyperoval of  $\cM_U(2^{2h})$.
\end{theorem}
We shall also see that  any oval arising from a  regular hyperoval
 of $AG(2,2^{2h})$ by deleting a point cannot be  inherited by
 $\cM_U(2^{2h})$.
The hypothesis on $\Omega$ being  a regular hyperoval
cannot be dropped; see \cite{PRS} for
examples of hyperovals in the Moulton plane of order $16$.




\par


\section{Proof of Theorem \ref{conf}}

We begin by showing the following lemma,
which is a slight generalisation of  Lemma 2.1. in \cite{AL}.
\begin{lemma}
\label{lea}
 Let $q$ be any prime power.
A pencil of affine lines $\cL(P)$ of $\cM_U(q^2)$ with centre
$P(x_0,y_0)$, either consists of lines of a Baer subplane $\cB$ of
$PG(2,q^2)$, or is a pencil in $AG(2,q^2)$ with the same center
according to $||x_0||\in U$ or not.  In particular, in the former
case, the $q^2+1$ lines in $\cL(P)$ plus the $q$ vertical lines $X=c$
with $||c||=x_0^{q+1}$ and $c\neq x_0$ are the lines of $\cB$.
\end{lemma}
\begin{proof}
The pencil $\cL(P)$ consists of the lines
\[r_m: y=m \odot x-m\odot x_0+y_0, \] with $m \in GF(q^2)$, plus the
vertical line $\ell: x=x_0$.  First suppose $||x_0||\in U$. In this
case $m\odot x_0= m^qx_0$ and the line $r_m$ of $\cL(P)$ corresponds
to the point $(m, m^qx_0-y_0)$ in the dual of $\cM_U(q^2)$, which is
an Andr\'e plane.

 As $m$ varies over $GF(q^2)$ we get $q^2$ affine points of the Baer
subplane $\cB'$ in $PG(2,q^2)$ represented by $y=x^qx_0-y_0$. The
points at infinity of $\cB'$ are those points $(c)$ such that
$c^{q+1}=||x_0||$.  As the dual of a Baer subplane is a Baer subplane,
it follows that the lines in $\cL(P)$ are the lines of a Baer subplane
$\cB$ in $PG(2,q^2)$. More precisely, the lines in $\cL(P)$ plus the
$q$ vertical lines $x=c$, $||c||=||x_0||$, with $c \neq x_0$, are the
lines of $\cB$.

In the case in which $||x_0||\notin U$ the line $r_m: y=m \odot
x-mx_0+y_0 $ in $\cL(P)$ corresponds to the point \[(m, mx_0-y_0)\] in
the dual of $\cM_U(q^2)$. As $m$ varies over $GF(q^2)$ we get $q^2$
affine points in $AG(2,q^2)$ on the line $y=x_0x-y_0$. The dual of
each such a point is a line of $AG(2,q^2)$ through the point
$P(x_0,y_0)$.
\end{proof}




Let $\Omega$  denote a  regular  hyperoval in  $AG(2,q^{2})$, $q=2^h$, $h>1$.
 It will be shown that for any point $P(x_0,y_0)$ with $||x_0||\in U$
there is at least a $2$--secant to $\Omega^{\star}=\Omega\cap \Phi$ in
$\cM_U(q^{2})$ through $P$.

 Assume $\cB$ to be the Baer subplane in $PG(2,q^{2})$ containing the
lines of $\cL(P)$ and the $q$ vertical lines $X=c$ with
$||c||=x_0^{q+1}$, $c\neq x_0$.  Write $\Delta$ for the set of all
points of $\Omega$ not covered by a vertical line of $\cB$ and also
let $n=|\Delta|$ and $m=q^2+2-n$.  The vertical lines of $\cB$ cover
at most $2(q+1)$ points of $\Omega$; thus, $q^2-2q\leq n\leq q^2+2$.
We shall show that there is at least a line in $\cB$ meeting $\Delta$
in two points.

\par
Let $T\in\Delta$; since $T\not\in\cB$, there is a unique line
$\ell_T$ of $\cB$ through $T$.
Every point $Q\in\Omega\setminus\Delta$  lies on at most
$q+1-(m-1)=q-m+2$ lines $\ell_T$ with $T\in\Delta$.
Suppose by contradiction that  for every $T\in\Delta$,
\[\ell_T \cap \Omega=\{ T, Q \}, \mbox{with} \ Q \in \Omega\setminus\Delta.\]
 The
total number of lines $\ell_T$ obtained as $Q$ varies in $\Omega\setminus\Delta$
does not exceed $m(q-m+2)$.
So,
\[ n=q^2-m+2\leq m(q-m+2). \]
As $m$ is a non--negative integer, this is possible only for $q=2$.

Since $\ell_T$ is not a vertical line, it turns out to be  a chord of $\Omega^*$ in $\cM_U(q^2)$  passing through $P(x_0,y_0)$.
This implies that no point $P(x_0,y_0) \in \Psi$ may be aggregated
to $\Omega^{\star}$ in order to obtain an arc.

This holds true in the case $P(x_0,y_0) \in \Phi$. In $AG(2,q^2)$
there pass $(q^2+2)/2$ secants to $\Omega$ through a point
$P(x_0,y_0)\notin \Omega$ and, hence, $N=(q^2+2)/2-s$ secants to
$\Omega^*$, where $s\leq 2(q+1)|U|$. So by the hypothesis $|U|<q/4-1$,
we obtain $N>0$; this implies that no point $P(x_0,y_0) \in \Phi$ may
be aggregated to $\Omega^{\star}$ in order to obtain a larger arc. The
same argument works also when $P$ is assumed to be a point at
infinity. Theorem \ref{conf} is thus proved.
\par

\section{Proof of Theorem \ref{maint}}
We shall use
the notion of
 \emph{conic blocking set}; see \cite{H}.
A conic blocking set $\cB$ is a set of lines in a
Desarguesian projective plane met by all conics;
a conic blocking
set $\cB$ is \emph{irreducible} if for any line of $\cB$ there is a
conic intersecting $\cB$ in just that line.

\begin{lemma}[Theorem 4.4,\cite{H}]
\label{lh} The line--set \[\cB=\{y=mx:m\in\GF(q)\}\cup\{x=0\}\] is
an irreducible conic blocking set in $\PG(2,q^2)$, where $q=2^{h}$, $h>1$.
\end{lemma}

\begin{lemma}
\label{main}
 Let $\Omega$ be a regular hyperoval of $PG(2,q^2)$, with  $q=2^h$, $h>1$. Then,
 there are at least two points $P(x,y)$ in $\Omega$ such that
 $||x||\in U$.
\end{lemma}
\begin{proof}
To prove the lemma we show that
the set $ \Psi'=\Psi \cup Y_{\infty}$, is
a conic blocking set.
We observe that the conic blocking set of Lemma \ref{lh}
is actually a degenerate Hermitian curve of $\PG(2,q^2)$ with
 equation $x^{q}y-xy^q=0$. Since all degenerate Hermitian curves
are projectively equivalent, this implies that any such a
curve is a conic blocking set.
On the other hand, $\Psi'$ may be regarded as the union of degenerate
Hermitian curves of equation $x^{q+1}=cz^{q+1}$, as $c$ varies in $U$.
Thus, $\Psi'$ is also a conic blocking set.
Suppose now $\Omega=\cC \cup N$, where $\cC$ is a conic of nucleus $N$.
Take $P\in \Psi'\cap\cC $. If $P=Y_{\infty}$ then at most one of the
vertical lines $X=c$, with $||c||\in U$, is tangent to $\cC$; hence
there are at least $q$ points $P'(x,y)\in \Omega $ with $||x||\in U$
and thus $|\Psi \cap \Omega|\geq q$.

Next, assume that $P=P(x,y) \in \Psi$. If the line $[x]$ is secant to
$\cC$ the assertion immediately follows.  If the line $[x]$ is tangent
to $\cC$ then the nucleus $N$ lies on $[x]$. Now, either $N$ is an
affine point in $\Psi$ or $N=Y_{\infty}$. In the former case we have
$|\Psi \cap \Omega|\geq 2$; in the latter, the lines $X=c$ with
$||c||\in U$ are all tangent to $\cC$; hence, there are at least other
$q+1$ points $P'(x,y)\in \Omega$ such that $||x||\in U$.

\end{proof}
\par Now, let $\Omega$ be a regular hyperoval in $AG(2,q^2)$, with
$q=2^h$ and $h>1$.  From Lemma \ref{main} we deduce that $|\Omega^*
\cap \Phi|\leq q^2 $; furthermore, Theorem \ref{maint} guaranties that
$\Omega^*$ is a complete arc in the projective closure of
$\cM_{U}(q^2)$, whence Theorem \ref{maint} follows.

\begin{remark} The largest arc of $\cM_U(q^2)$ contained in a regular
hyperoval of $AG(2,q^2)$, with $q=2^h$, has at most $q^2$ points; in
particular any oval which arises from a hyperoval of $AG(2,q^2)$ by
deleting a point cannot be an oval of $\cM_U(q^2)$.  For an actual
example of a $q^2$--arc of $\cM_U(q^2)$ coming from a regular
hyperoval of $AG(2,q^2)$ see \cite{M}.  This also shows that the
result of \cite{kg1974} cannot be extended to even $q$.
\end{remark}

\medskip
\begin{minipage}{7cm}
{\sc Angela Aguglia}\\
Dipartimento di Matematica \\
Politecnico di Bari \\
Via Orabona 4 \\
I--70125 Bari \\
Italy \\
e-mail: {\begin{minipage}[t]{4cm}
    {\tt a.aguglia@poliba.it}
        \end{minipage}}
\end{minipage}
\vskip1cm\noindent
\begin{minipage}{7cm}
{\sc Luca Giuzzi}\\
Dipartimento di Matematica \\
Universit\`a degli Studi di Brescia\\
Facolt\`a di Ingegneria \\
Via Valotti 9 \\
I-25133 Brescia \\
Italy \\
e-mail: {\begin{minipage}[t]{4cm}
            {\tt giuzzi@ing.unibs.it}
        \end{minipage}}
\end{minipage}

\begin{thebibliography}{999}

\bibitem{AL} V. Abatangelo, B. Larato, \emph{Canonically Inherited
    Arcs in Moulton Planes of Odd Order}, to appear on Innovations in Incidence Geometry.

\bibitem{HCD} C.J. Colbourn, J.H. Diniz, \emph{Handbook of
    Combinatorial Designs}, Chapman \& Hall (2007).
\bibitem{D} P. Dembowski, \emph{Finite Geometries}, Springer (1968).
\bibitem{JJB} N.L. Johnson, V. Jha, M. Biliotti, {\em Handbook of finite translation planes}. Pure and Applied Mathematics (Boca Raton), {\bf 289}
 Chapman \& Hall/CRC, Boca Raton, FL, (2007).
\bibitem{kg1974} G. Korchm\'aros,
  {\em Ovali nei piani di Moulton d'ordine dispari},
  Colloquio Internazionale sulle Teorie Combinatorie (Rome, 1973),
  Tomo II,  395--398. Atti dei Convegni Lincei,
  No. {\bf 17}, Accad. Naz. Lincei, Rome,  (1976).
\bibitem{kg1976} G. Korchm\'aros, \emph{Ovali nei piani di Hall di
    ordine dispari},
 Atti Accad. Naz. Lincei Rend. Cl. Sci. Fis. Mat. Natur. (8)
{\bf 56} (1974), no. 3, 315--317.

\bibitem{H} L. D. Holder,
  {\em Conic blocking sets in Desarguesian projective planes},
  J. Geom. {\bf 80} (2004), no. {1-2}, 95--105.

\bibitem{HP} D.Hughes, F. Piper, \emph{Projective planes}, Springer (1973).

\bibitem{M} G. Menichetti,
  {\em $k$--archi completi in piani di Moulton d'ordine $q^m$},
  Atti Accad. Naz. Lincei Rend. Cl. Sci. Fis. Mat. Natur. (8) {\bf 60}
  (1976),  no. {6}, 775--781.
\bibitem{Mo} F.R. Moulton, \emph{A simple non--desarguesian plane
    geometry}, Trans. Am. Math. Soc. {\bf 3} (1902), 192--195.

\bibitem{PRS} T. Penttila, G.F. Royle, M.K. Simpson,
  \emph{Hyperovals in the known projective planes of order $16$},
  J. Combin. Des. {\bf 4} (1996), no. {1}, 59--65.
\end{thebibliography}
\end{document}